\documentclass[11pt,a4paper,leqno]{article}
\usepackage{graphicx}
\usepackage{amsmath,amssymb}
\usepackage{amsthm}
\usepackage{newpxtext,newpxmath}
\title {A newly-generalized problem from a problem for the Mathematical Olympiad and the methods to solve it. }
\author {Yasushi Ieno}
\date{}
\begin{document}
\maketitle 
% \newpage
\noindent Abstract. A newly-generalized problem from a problem initially thought for the Mathematical Olympiad and the methods to solve it. \\\\
% 2020 Mathematics Subject Classification. Primary 05A05; Secondary 05A18, 11B37.\\\\
\\\\
% \noindent AMS classification.  05A05, 05A20.\\
% \begin{figure}[htbp]
% \begin{center}
% \includegraphics [width=80mm]{CIMG1137.JPG}
% \caption{Author}
% \end{center}
% \end{figure}
% $\{a_1,a_2,\cdots,a_n\}$\\
% \includegraphics[bb={0 0 595 842}]{maxmin.pdf}
0. Introduction\\

Recently I read [1], in which Dvornicich, Veneziano and Zannier show a problem initially thought for the Romanian Mathematical Olympiad 2010 (see [2]). They present its two solutions, with the aid of its several interpretations. The one is an orthodox method, which is essentially the same as what the author of this problem (see [1]) has shown. And the other is a method in which they generalize the recurrence sequences in this problem to recurrence sequences related to a much longer set of diophantine equations, and use a linear differential equation of the second order. 

I have generalized this problem into a new problem and searched for its solutions by the reference of [1].
\\\\
1. An original problem shortlisted for the Mathematical Olympiad and its newly-generalized problem.\\

\noindent \ \ 1-1 An original problem shortlisted for the Mathematical Olympiad.\\

{\bf Problem 1} \ Let $\rm x_0$, $\rm x_1$, $\rm x_2$, . . . be the sequence defined by
\begin{eqnarray}
\lefteqn {\rm x_0\ {=}\ 1} \ \ \ \ \ \ \ \ \ \ \ \ \ \ \ \ \ \ \ \ \ \ \ \ \ \ \ \ \ \ \ \ \ \ \ \ \ \ \ \ \nonumber \\
\lefteqn {\rm x_{n+1}\ {=}\ 1\ {+}\ \frac {n}{x_n} \ \ \ \forall n \geq 0.} \ \ \ \ \ \ \ \ \ \ \ \ \ \ \ \ \ \ \ \ \ \ \ \ \ \ \ \ \ \ \ \ \ \ \ \ \ \ \ \ 
\end{eqnarray}

What are the values of n for which ${\rm x_n}$ is an integer?
\\\\

\noindent \ \ 1-2 A generalized problem from a problem for the Mathematical Olympiad.\\

{\bf Problem 2} \ Let $\rm x_0$, $\rm x_1$, $\rm x_2$, . . . be the sequence defined by
\begin{eqnarray}
\lefteqn {\rm x_0\ {=}\ 1}\ \ \ \ \ \ \ \ \ \ \ \ \ \ \ \ \ \ \ \ \ \ \ \ \ \ \ \ \ \ \ \ \ \ \ \ \ \ \ \ \ \ \ \ \ \ \ \ \ \ \ \ \ \ \ \ \ \nonumber \\
\lefteqn {\rm x_{n+1}\ {=}\ c\ {+}\ \frac {n}{x_n} \ \forall n \geq 0, where\ c\ is\ a\ constant\ positive\ integer.}\ \ \ \ \ \ \ \ \ \ \ \ \ \ \ \ \ \ \ \ \ \ \ \ \ \ \ \ \ \ \ \ \ \ \ \ \ \ \ \ \ \ \ \ \ \ \ \ \ \ \ \ \ \ \ \ \ \ \ \ 
\end{eqnarray}

What are the values of n for which ${\rm x_n}$ is an integer?\\\\

As mentioned above, we have generalized (1) into (2).

From now on we will research how to solve {\bf Problem 2} by the reference of [1].\\

\noindent 2. <First Solution>.
% Note that 
% \newpage 
\newtheorem{theorem}{Theorem} % theorem 環境を定義
% \begin{theorem} % Theorem 1 
% \begin{eqnarray}
% \lefteqn {{\rm  If \ n \geq 2 \ then}} \nonumber \\
% && \left| {\rm C(n)} \right|{=}(2^{\rm n-1}{-}1)\left| {\rm P_C(n,1)} \right| {+} \sum_{\rm k{=}2}^{\rm n-1} ({-}(1/2)^{\rm k-1}\left| {\rm % P(n,k)} \right|{+}2^{\rm n-k}\left| {\rm P_C(n,k)} \right|){+}\left| {\rm P_C(n,n)} \right| 
% \end{eqnarray}
% \end{theorem}
\\

Now look at (2). 

If we define ${\rm f_n(x)}$=c+${\rm \frac{x}{n}}$, we see that ${\rm x_{\rm n{+}1}}$=${\rm f_n(x_n)}$, so that one is led to
study the dynamics of the sequence of functions  ${\rm f_n}$.

Let us call ${\rm y_n}$ = ${\rm \frac{c+\sqrt{4n+c^2}}{2}}$ the (positive) fixed point of ${\rm f_n}$. We have that if x<${\rm y_n}$ then ${\rm f_n(x)}$>${\rm y_n}$ and vice versa.

From now on we newly add the condition for (2) that c $\geq$ 2, because the solutions for c=1 are already said in [1].  

We can prove by induction that\\

{\bf Lemma 1} For every n $\geq$ 4 we have
\begin{eqnarray}
\lefteqn {\rm y_{n{-}1}{=}\frac{c{+}\sqrt{4(n{-}1){+}c^2}}{2}<x_n<\frac{c{+}\sqrt{4n{+}c^2}}{2}{=}y_n} \ \ \ \ \ \ \ \ \ \ \ \ \ \ \ \ \ \ \ \ \ \ \ \ \ \ \ \ \ \ \ \ \ \ \ \ \ \ \ \ \ \ \ \ \ \ \ \ \ \ \ \ \ \ \ \ \ \ \ \ \ \ \ \ \ \ \ \ \ \ \ \ \ \ \ 
\end{eqnarray}
\begin{proof}
% \quad\par
By a direct computation, we have ${\rm \frac{c{+}\sqrt{4{+}c^2}}{2}}$<c<${\rm \frac{c{+}\sqrt{8{+}c^2}}{2}}$, which establishes
the basis of the induction. Assuming (3) holds for n, by the previous remark we have that ${\rm y_n}$ < ${\rm x_{\rm n+1}}$, so we need only to prove that
\begin{eqnarray}
\lefteqn {\rm c{+}\frac{n}{x}<y_{\rm n+1}} \ \ \ \ \ \ \ \ \ \ \ \ \ \ \ \ \ \ \ \ \ \nonumber
\end{eqnarray}

By the inductive hypothesis, it is enough to show that
\begin{eqnarray}
\lefteqn {\rm c{+}\frac{n}{y_{\rm n{-}1}}<y_{\rm n{+}1}\  i.e.\ ,} \ \ \ \ \ \ \ \ \ \ \ \ \ \ \ \ \ \ \ \ \ \ \nonumber \\
\lefteqn {\rm c{+}\frac{2n}{c{+}\sqrt{4n{+}(c^2{-}4)}}<\frac{c{+}\sqrt{4n{+}(c^2{+}4)}}{2},} \ \ \ \ \ \ \ \ \ \ \ \ \ \ \ \ \ \ \ \ \ \  \nonumber
\end{eqnarray}
which is an elementary, though tedious, computation. To show it we had better use 
\begin{eqnarray}
\lefteqn  {\rm \sqrt{1{+}x}{\geq}{1{+}\frac{x}{2}}\ \ \ \forall x{>}0.} \ \ \ \ \ \ \ \ \ \ \ \ \ \ \ \ \ \ \ \ \ \ \ \ \ \ \ \ \ \ \ \ \ \ \ \ \ \ \ \ \ \ \ \ \ \ \ \ \ \ \ \ \ \ \ \ \ \ \ \ \ \ \ \ \ \ \ \ \ \ \ \ \nonumber
\end{eqnarray}
\end{proof}
\noindent Remark 1. For the values of n smaller than 2, we have ${\rm y_0}$=${\rm x_1}$=c and ${\rm y_1}$=${\rm \frac{c{+}\sqrt{c^2{+}4}}{2}}$.\\

Let us now assume that ${\rm x_n}$ is an integer for some n $\geq$ 4. From the lemma we have
\begin{eqnarray}
\lefteqn {\rm \frac{c{+}\sqrt{4(n{-}1){+}c^2}}{2}<x_n<\frac{c{+}\sqrt{4n{+}c^2}}{2}}  \ \ \ \ \ \ \ \ \ \ \ \ \ \ \ \ \ \ \ \ \ \ \ \ \ \ \ \ \ \ \ \ \ \ \ \ \ \ \ \ \ \ \ \ \ \ \ \ \ \ \ \ \ \ \ \ \ \ \ \ \ \ \ \ \ \ \ \ \ \ \ \ \ \ \ \nonumber \\
\lefteqn {\rm \sqrt{4(n{-}1){+}c^2}<2x_n{-}c<\sqrt{4n{+}c^2}}  \ \ \ \ \ \ \ \ \ \ \ \ \ \ \ \ \ \ \ \ \ \ \ \ \ \ \ \ \ \ \ \ \ \ \ \ \ \ \ \ \ \ \ \ \ \ \ \ \ \ \ \ \ \ \ \ \ \ \ \ \ \ \ \ \ \ \ \ \ \ \ \ \ \ \ \nonumber \\
\lefteqn {\rm {4(n{-}1){+}c^2}<(2x_n{-}c)^2<{4n{+}c^2}}  \ \ \ \ \ \ \ \ \ \ \ \ \ \ \ \ \ \ \ \ \ \ \ \ \ \ \ \ \ \ \ \ \ \ \ \ \ \ \ \ \ \ \ \ \ \ \ \ \ \ \ \ \ \ \ \ \ \ \ \ \ \ \ \ \ \ \ \ \ \ \ \ \ \ \ \nonumber
\end{eqnarray}

However the last inequalities are inconsistent modulo 4, whether c is even or odd.

Therefore we conclude that the only integral values of the sequence are ${\rm x_0}$, ${\rm x_1}$.\\

Essentially the same solution may be reached by a slightly different approach.

The same conclusion as before can be reached if we show that ${\rm n{-}1}$<${\rm x_n^2{-}cx_n}$<n for n$\geq$2.

We argue by induction. The inequalities are verified by direct inspection for n${=}$2, since\\
 1<${\rm {(c{+}\frac{1}{c})}^2{-}c(c{+}\frac{1}{c}){=}1{+}\frac{1}{c^2}}$<2.\\

Now let an = ${\rm x_n^2-cx_n}$, and assume that the inequalities hold up to n. We may write ${\rm a_{n-1}}$ as
\begin{eqnarray}
\lefteqn {\rm a_{n+1}{=}x_{n+1}^2{-}cx_{n+1}{=}x_{n+1}(x_{n+1}{-}c){=}(c{+}\frac{n}{x_n})\frac{n}{x_n}{=}\frac{n(cx_n{+}n)}{x_n^2}} \ \ \ \ \ \ \ \ \ \ \ \ \ \ \ \ \ \ \ \ \ \ \ \ \ \ \ \ \ \ \ \ \ \ \ \ \ \ \ \ \ \ \ \ \ \ \ \ \ \ \ \ \ \ \ \ \ \ \ \ \ \ \ \ \ \ \ \ \ \ \ \ \ \ \ \nonumber
\end{eqnarray}

By the induction hypothesis we have:
\begin{eqnarray}
\lefteqn {\rm x_n^2{<}cx_n{+}n\Rightarrow a_{n{+}1}{>}n\frac{cx_n{+}n}{cx_n{+}n}{=}n}  \ \ \ \ \ \ \ \ \ \ \ \ \ \ \ \ \ \ \ \ \ \ \ \ \ \ \ \ \ \ \ \ \ \ \ \ \ \ \ \ \ \ \ \ \ \ \ \ \ \ \ \ \ \ \ \ \ \ \ \ \ \ \ \ \ \ \ \ \ \ \ \ \ \ \ \nonumber 
\end{eqnarray}
and
\begin{eqnarray}
\lefteqn {\rm x_n^2{>}cx_n{+}n{-}1\Rightarrow a_{n{+}1}{<}n\frac{cx_n{+}n}{cx_n{+}n{-}1}{<}n{+}1}  \ \ \ \ \ \ \ \ \ \ \ \ \ \ \ \ \ \ \ \ \ \ \ \ \ \ \ \ \ \ \ \ \ \ \ \ \ \ \ \ \ \ \ \ \ \ \ \ \ \ \ \ \ \ \ \ \ \ \ \ \ \ \ \ \ \ \ \ \ \ \ \ \ \ \ \nonumber
\end{eqnarray}
since ${\rm cx_n{>}1}$.\\

\noindent 3. <Second Solution>.\\

To study the sequence (${\rm x_n}$) from an arithmetic point of view we define an integer sequences (${\rm a_n}$) by the recurrences
\begin{eqnarray}
\lefteqn {\rm a_0{=}1}  \ \ \ \ \ \ \ \ \ \ \ \ \ \ \ \ \ \ \ \ \ \ \ \ \ \ \ \ \ \ \ \ \ \ \ \ \ \ \ \ \ \ \ \ \ \ \ \ \ \ \ \ \ \ \ \ \ \ \ \ \ \ \ \ \ \ \ \ \ \ \ \ \ \ \ \nonumber \\
\lefteqn {\rm a_1{=}c}  \ \ \ \ \ \ \ \ \ \ \ \ \ \ \ \ \ \ \ \ \ \ \ \ \ \ \ \ \ \ \ \ \ \ \ \ \ \ \ \ \ \ \ \ \ \ \ \ \ \ \ \ \ \ \ \ \ \ \ \ \ \ \ \ \ \ \ \ \ \ \ \ \ \ \ \nonumber \\
\lefteqn {\rm a_{n{+}2}{=}ca_{n{+}1}{+}(n{+}1)a_n \ \ \ \forall n \geq 0, }\ \ \ \ \ \ \ \ \ \ \ \ \ \ \ \ \ \ \ \ \ \ \ \ \ \ \ \ \ \ \ \ \ \ \ \ \ \ \ \ \ \ \ \ \ \ \ \ \ \ \ \ \ \ \ \ \ \ \ \ \ \ \ \ \ \ \ \ \ \ \ \ \ \ \  
\end{eqnarray}
with ${\rm x_n{=}\frac{a_n}{a_{n{-}1}}}$.\\

Let us define ${\rm d_n{=}gcd(a_n,a_{n{-}1});\ d_n}$ tells us how much the reduced denominator of ${\rm x_n}$ differs from ${\rm a_{n{-}1}}$. So, to obtain a lower bound for said denominator, we need a lower bound for ${\rm a_{n{-}1}}$ and an upper bound for ${\rm d_n}$.
\\\\
Remark 2. By the recurrence (4) we see that ${\rm d_{n{+}1}}$|${\rm a_{n{+}2}}$, and so ${\rm d_{n{+}1}}$|${\rm d_{n{+}2}}$; this will be helpful in establishing an upper bound for ${\rm d_n}$.\\

A lower bound for ${\rm a_n}$ is easily obtained as in the following lemma.
\\\\
{\bf Lemma 2} For every n$\geq$0 we have ${\rm a_n}$$\geq$$\sqrt{n!}$.

\begin{proof} We argue by induction on n$\geq$0. We check that ${\rm a_0}$=1=$\sqrt{0!}$, ${\rm a_1}$=c$\geq$$\sqrt{1!}$, and assuming the bound for ${\rm a_n}$ and ${\rm a_{n{+}1}}$ we get
\begin{eqnarray}
\lefteqn {\rm a_{n{+}2}{=}ca_{n{+}1}{+}(n{+}1)a_n\geq \sqrt{c(n{+}1)!}{+}(n{+}1)\sqrt{n!}}  \ \ \ \ \ \ \ \ \ \ \ \ \ \ \ \ \ \ \ \ \ \ \ \ \ \ \ \ \ \ \ \ \ \ \ \ \ \ \ \ \ \ \ \ \ \ \ \ \ \ \ \ \ \ \ \ \ \ \ \ \ \ \ \ \ \ \ \ \ \ \ \ \ \ \ \nonumber \\
\lefteqn {\rm {=}\sqrt{(n{+}2)!}\ \frac{c{+}\sqrt{n{+}1}}{\sqrt{n{+}2}}{=}\sqrt{(n{+}2)!}\sqrt{\frac{n{+}1{+}c^2}{n{+}2}{+}\frac{2c\sqrt{n{+}1}}{n{+}2}}}\ \ \ \ \ \ \ \ \ \ \ \ \ \ \ \ \ \ \ \ \ \ \ \ \ \ \ \ \ \ \ \ \ \ \ \ \ \ \ \ \ \ \ \ \ \ \ \ \ \ \ \ \ \ \ \ \ \ \ \ \ \ \ \ \ \ \ \ \ \ \ \ \ \ \ \ \ \ \nonumber \\
\lefteqn {\rm \geq \sqrt{(n{+}2)!}}.  \ \ \ \ \ \ \ \ \ \ \ \ \ \ \ \ \ \ \ \ \ \ \ \ \ \ \ \ \ \ \ \ \ \ \ \ \ \ \ \ \ \ \ \ \ \ \ \ \ \ \ \ \ \ \ \ \ \ \ \ \ \ \ \ \ \ \ \ \ \ \ \ \ \ \ \nonumber 
\end{eqnarray}
\end{proof}
To get an upper bound for ${\rm d_n}$ we introduce the exponential generating function of the sequence ${\rm(a_n)_{n\in\mathbb{N}}}$, namely
\begin{eqnarray}
\lefteqn {\rm F(x){=}\sum_{n=0}^{\infty}\frac{a_n}{n!}x^n} \ \ \ \ \ \ \ \ \ \ \ \ \ \ \ \ \ \ \ \ \ \ \ \ \ \ \ \ \ \ \ \ \ \ \ \ \ \ \ \ \ \ \ \ \ \ \ \ \ \ \ \ \ \ \ \ \ \ \ \ \ \ \ \ \ \ \ \ \ \ \ \ \ \ \ \nonumber 
\end{eqnarray}

We consider F(x) merely as a formal power series, although one could prove that it converges for every complex x. 
From the recurrence on (${\rm a_n}$) we can obtain a differential equation for F; in fact, we can multiply (6) by ${\rm \frac{a_n}{n!}}$ and sum it for n$\geq$0; since clearly F'(x)=${\rm \sum_{n=0}^{\infty}\frac{a_{n{+}1}}{n!}x^n}$ and F''(x)=${\rm \sum_{n=0}^{\infty}\frac{a_{n{+}2}}{n!}x^n}$, we obtain that F satisfies the conditions
\begin{eqnarray}
\left\{ \begin{array}{l}
{\rm F(0){=}1} \\
{\rm F'(0){=}c} \\
{\rm F''(x){=}(c{+}x)F'(x){+}F(x)}
\end{array} \right.
\end{eqnarray}

The Cauchy problem (5) may be solved (in the ring of formal power series) to get
\begin{eqnarray}
\lefteqn {\rm F(x){=}e^{cx{+}\frac{x^2}{2}}}  \ \ \ \ \ \ \ \ \ \ \ \ \ \ \ \ \ \ \ \ \ \ \ \ \ \ \ \ \ \ \ \ \ \ \ \ \ \ \ \ \ \ \ \ \ \ \ \ \ \ \ \ \ \ \ \ \ \ \ \ \ \ \ \ \ \ \ \ \ \ \ \ \ \ \ \nonumber
\end{eqnarray}
and we can use this explicit form to get a formula for ${\rm a_n}$. In fact
\begin{eqnarray}
\lefteqn {\rm \sum_{n=0}^{\infty}\frac{a_n}{n!}x^n{=}e^{cx{+}\frac{x^2}{2}}{=}e^{cx}e^{\frac{x^2}{2}}}  \ \ \ \ \ \ \ \ \ \ \ \ \ \ \ \ \ \ \ \ \ \ \ \ \ \ \ \ \ \ \ \ \ \ \ \ \ \ \ \ \ \ \ \ \ \ \ \ \ \ \ \ \ \ \ \ \ \ \ \ \ \ \ \ \ \ \ \ \ \ \ \ \ \ \nonumber \\
\lefteqn {\rm {=}\sum_{m=0}^{\infty}\frac{{c^m}{x^m}}{m!}\sum_{s=0}^{\infty}\frac{x^{2s}}{{2^s}s!}}  \ \ \ \ \ \ \ \ \ \ \ \ \ \ \ \ \ \ \ \ \ \ \ \ \ \ \ \ \ \ \ \ \ \ \ \ \ \ \ \ \ \ \ \ \ \ \ \ \ \ \ \ \ \ \ \ \ \ \ \ \ \ \ \ \ \ \ \ \ \ \ \ \ \ \nonumber \\
\lefteqn {\rm \frac{a_n}{n!}{=}\sum_{2s{+}m=n}^{}\frac{c^m}{{2^s}m!s!}}  \ \ \ \ \ \ \ \ \ \ \ \ \ \ \ \ \ \ \ \ \ \ \ \ \ \ \ \ \ \ \ \ \ \ \ \ \ \ \ \ \ \ \ \ \ \ \ \ \ \ \ \ \ \ \ \ \ \ \ \ \ \ \ \ \ \ \ \ \ \ \ \ \ \ \nonumber \\
\lefteqn {\rm a_n{=}\sum_{2s\leq n}^{}\frac{n!}{{c^{2s{-}n}}{2^s}(n{-}2s)!s!}}  \ \ \ \ \ \ \ \ \ \ \ \ \ \ \ \ \ \ \ \ \ \ \ \ \ \ \ \ \ \ \ \ \ \ \ \ \ \ \ \ \ \ \ \ \ \ \ \ \ \ \ \ \ \ \ \ \ \ \ \ \ \ \ \ \ \ \ \ \ \ \ \ \ \ \nonumber \\
\lefteqn {\rm {=}\sum_{2s\leq n}^{}c^{n{-}2s}{\binom{n}{2s}}(2s{-}1)!!}  \ \ \ \ \ \ \ \ \ \ \ \ \ \ \ \ \ \ \ \ \ \ \ \ \ \ \ \ \ \ \ \ \ \ \ \ \ \ \ \ \ \ \ \ \ \ \ \ \ \ \ \ \ \ \ \ \ \ \ \ \ \ \ \ \ \ \ \ \ \ \ \ \ \ \nonumber 
% \lefteqn {\rm {=}\sum_{m=0}^{\infty}\frac{{x^m}{m!}}} \ \ \ \ \ \nonumber 
\end{eqnarray}
where the semifactorial (${\rm 2s{-}1}$)!! denotes as usual the product (${\rm 2s{-}1}$)·(${\rm 2s{-}3}$)· · · 3·1 and is defined to be 1 for s=0.

We can now use the preceding formula to prove the following lemma.
\\\\
{\bf Lemma 3.} Let p be an odd prime. If p$\mid$n, then ${\rm a_n\equiv c^n}$(mod p).
\begin{proof}
Let p be an odd prime dividing n.\\
If p<2s, we have that p${\mid}$(2s${-}$1)!!, as p itself is one of the factors in the defining product of (2s${-}$1)!!.\\
If 0<2s<p the binomial ${\rm \binom{n}{2s}}$ is divisible by p, as the p factor in n is not cancelled by (2s)!.\\
So we have that in formula (8) only the term with s = 0 is not divisible by p, whence
\begin{eqnarray}
\lefteqn {\rm a_n{=}\sum_{2s\leq n}^{}c^{n{-}2s}{\binom{n}{2s}}(2s{-}1)!!{\equiv}{c^n}(mod\ p)}  \ \ \ \ \ \ \ \ \ \ \ \ \ \ \ \ \ \ \ \ \ \ \ \ \ \ \ \ \ \ \ \ \ \ \ \ \ \ \ \ \ \ \ \ \ \ \ \ \ \ \ \ \ \ \ \ \ \ \ \ \ \ \ \ \ \ \ \ \ \ \ \ \ \ \nonumber 
\end{eqnarray}
\end{proof}
Applying this lemma we get the following property of ${\rm d_n}$,
\\\\
{\bf Corollary 4.} For every n$\geq$1, ${\rm d_n}$ is a power of 2 and all of the odd prime factors of c.
\begin{proof}
If an odd prime p, that is not a prime factor of c, divides ${\rm d_m}$ for some m$\geq$1. Then, by Remark 2, p divides
${\rm d_n}$ (and hence ${\rm a_n}$) for all n$\geq$m, so also for n=pm. But this is not possible because ${\rm a_pm}$$\equiv$1(mod p) by Lemma 3.
\end{proof}

We are now ready to prove an upper bound for ${\rm d_n}$, which will follow by using again the exponential generating function F(x).
\\\\
{\bf Proposition 5.} For every n$\geq$2 we have that 
\begin{eqnarray}
\lefteqn {\rm d_n\leq2^{n{-}1}\ (p_1^{\frac{1}{{p_1}{-}1}}p_2^{\frac{1}{{p_2}{-}1}}. . .p_j^{\frac{1}{{p_j}{-}1}})^{n{-}2}\  ({2n{-}3})^{j/2}} \ \ \ \ \ \ \ \ \ \ \ \ \ \ \ \ \ \ \ \ \ \ \ \ \ \ \ \ \ \ \ \ \ \ \ \ \ \ \ \ \ \ \ \ \ \ \ \ \ \ \ \ \ \ \ \ \ \ \ \ \ \ \ \ \ \ \ \ \ \ \ \ \ \ \ \nonumber
\end{eqnarray}

where ${\rm p_1}$,${\rm p_2}$,. . .,${\rm p_j}$ are all of the distinct odd prime factors of c.  
\begin{proof}
We have the following identities concerning the above generating function F(x):
\begin{eqnarray}
\lefteqn {\rm (\sum_{m=0}^{\infty}a_m\frac{x^m}{m!})(\sum_{r=0}^{\infty}(-1)^r a_r\frac{x^r}{r!}){=}F(x)F(-x){=}e^{(x^2)}{=}\sum_{n=0}^{\infty}\frac{2^n}{n!}}  \ \ \ \ \ \ \ \ \ \ \ \ \ \ \ \ \ \ \ \ \ \ \ \ \ \ \ \ \ \ \ \ \ \ \ \ \ \ \ \ \ \ \ \ \ \ \ \ \ \ \ \ \ \ \ \ \ \ \ \ \ \ \ \ \ \ \ \ \ \ \ \ \ \ \nonumber 
\end{eqnarray}

Comparing the coefficients of ${\rm x_{2n}}$ for any n$\geq$1, we obtain
\begin{eqnarray}
\lefteqn {\rm \sum_{m{+}r=2n}^{}(-1)^r{\binom{2n}{m}}a_m a_r{=}\frac{(2n)!}{n!}{=}2^n(2n{-}1)!!}  \ \ \ \ \ \ \ \ \ \ \ \ \ \ \ \ \ \ \ \ \ \ \ \ \ \ \ \ \ \ \ \ \ \ \ \ \ \ \ \ \ \ \ \ \ \ \ \ \ \ \ \ \ \ \ \ \ \ \ \ \ \ \ \ \ \ \ \ \ \ \ \ \ \  
\end{eqnarray}

Now, as observed in Remark 3 above, ${\rm d_{n{+}1}}$ divides any ${\rm a_m}$ with m$\geq$n, so it divides the left-hand side of (6), and we know from the Corollary 4 that it is a power of 2 and all of the odd prime factors of c.

Therefore at first, the exponent of 2 in ${\rm d_{n{+}1}}$ is less than or equal to n.
Next, as mentioned above, we assume that c has as many as j distinct odd prime factors, ${\rm p_1}$,${\rm p_2}$,. . .,${\rm p_j}$.

Now we define ${\rm m_i}$ as the upper bound of the exponent of ${\rm p_i}$ in the prime factorization of ${\rm d_{n{+}1}}$, such that 1$\leq$i$\leq$j, then because ${\rm (2n{-}1)!!}$ is odd,
\begin{eqnarray}
\lefteqn {\rm m_i{=}\sum_{k{=}1}^{t_i}{\left[{\frac{(2n{-}1){+}p_i^k}{2p_i^k}}\right]}}  \ \ \ \ \ \ \ \ \ \ \ \ \ \ \ \ \ \ \ \ \ \ \ \ \ \ \ \ \ \ \ \ \ \ \ \ \ \ \ \ \ \ \ \ \ \ \ \ \ \ \ \ \ \ \ \ \ \ \ \ \ \ \ \ \ \ \ \ \ \ \ \ \ \  
\end{eqnarray}
where ${\rm t_i}${=}max t$(\in{\mathbb{R}})$ such that
\begin{eqnarray}
\lefteqn {\rm {\left[{\frac{(2n{-}1){+}p_i^t}{2p_i^t}}\right]}{\geq}1}. \ \ \ \ \ \ \ \ \ \ \ \ \ \ \ \ \ \ \ \ \ \ \ \ \ \ \ \ \ \ \ \ \ \ \ \ \ \ \ \ \ \ \ \ \ \ \ \ \ \ \ \ \ \ \ \ \ \ \ \ \ \ \ \ \ \ \ \ \ \ \ \ \ \  
\end{eqnarray}

It follows by (7) and (8) that
\begin{eqnarray}
\lefteqn {\rm m_i{=}\sum_{k{=}1}^{t_i}{\left[{\frac{2n{-}1}{2p_i^k}{+}\frac{1}{2}}\right]}}  \ \ \ \ \ \ \ \ \ \ \ \ \ \ \ \ \ \ \ \ \ \ \ \ \ \ \ \ \ \ \ \ \ \ \ \ \ \ \ \ \ \ \ \ \ \ \ \ \ \ \ \ \ \ \ \ \ \ \ \ \ \ \ \ \ \ \ \ \ \ \ \ \ \ \ \ \ \ \nonumber \\
\lefteqn {\rm \leq\sum_{k{=}1}^{t_i}{{\left(\frac{2n{-}1}{2p_i^k}{+}\frac{1}{2}\right)}}{=}\frac{2n{-}1}{2(p_i{-}1)}\left(1{-}{\frac{1}{{p_i}^{t_i}}}\right){+}\frac{t_i}{2} } \ \ \ \ \ \ \ \ \ \ \ \ \ \ \ \ \ \ \ \ \ \ \ \ \ \ \ \ \ \ \ \ \ \ \ \ \ \ \ \ \ \ \ \ \ \ \ \ \ \ \ \ \ \ \ \ \ \ \ \ \ \ \ \ \ \ \ \ \ \ \ \ \ \  \nonumber \\
\lefteqn {\rm \leq\frac{2n{-}1}{2(p_i{-}1)}\left(1{-}{\frac{1}{2n{-}1}}\right){+}\frac{t_i}{2}{=}\frac{n{-}1}{p_i{-}1}{+}\frac{t_i}{2}}  \ \ \ \ \ \ \ \ \ \ \ \ \ \ \ \ \ \ \ \ \ \ \ \ \ \ \ \ \ \ \ \ \ \ \ \ \ \ \ \ \ \ \ \ \ \ \ \ \ \ \ \ \ \ \ \ \ \ \ \ \ \ \ \ \ \ \ \ \ \ \ \ \ \  \nonumber \\
\lefteqn {\rm \therefore d_{n{+}1}\leq2^n p_1^{m_1}p_2^{m_2}. . .p_{j}^{m_j}}  \ \ \ \ \ \ \ \ \ \ \ \ \ \ \ \ \ \ \ \ \ \ \ \ \ \ \ \ \ \ \ \ \ \ \ \ \ \ \ \ \ \ \ \ \ \ \ \ \ \ \ \ \ \ \ \ \ \ \ \ \ \ \ \ \ \ \ \ \ \ \ \ \ \ \ \ \ \ \ \ \ \ \ \nonumber \\
\lefteqn {\rm \leq2^n\ p_1^{\frac{n{-}1}{{p_1}{-}1}{+}\frac{t_1}{2}}p_2^{\frac{n{-}1}{{p_2}{-}1}{+}\frac{t_2}{2}}. . .p_j^{\frac{n{-}1}{{p_j}{-}1}{+}\frac{t_j}{2}}}  \ \ \ \ \ \ \ \ \ \ \ \ \ \ \ \ \ \ \ \ \ \ \ \ \ \ \ \ \ \ \ \ \ \ \ \ \ \ \ \ \ \ \ \ \ \ \ \ \ \ \ \ \ \ \ \ \ \ \ \ \ \ \ \ \ \ \ \ \ \ \ \ \ \nonumber \\
\lefteqn {\rm {=}2^n\ (p_1^{\frac{1}{{p_1}{-}1}}p_2^{\frac{1}{{p_2}{-}1}}. . .p_j^{\frac{1}{{p_j}{-}1}})^{n{-}1}\  ({p_1^\frac{t_1}{2}}{p_2^\frac{t_2}{2}}. . .{p_j^\frac{t_j}{2}}) }  \ \ \ \ \ \ \ \ \ \ \ \ \ \ \ \ \ \ \ \ \ \ \ \ \ \ \ \ \ \ \ \ \ \ \ \ \ \ \ \ \ \ \ \ \ \ \ \ \ \ \ \ \ \ \ \ \ \ \ \ \ \ \ \ \ \ \ \ \ \ \ \ \ \nonumber \\
\lefteqn {\rm \leq2^n\ (p_1^{\frac{1}{{p_1}{-}1}}p_2^{\frac{1}{{p_2}{-}1}}. . .p_j^{\frac{1}{{p_j}{-}1}})^{n{-}1}\  ({2n{-}1})^{j/2} \ \ \ ({\rm \ see\ (8)\ })}  \ \ \ \ \ \ \ \ \ \ \ \ \ \ \ \ \ \ \ \ \ \ \ \ \ \ \ \ \ \ \ \ \ \ \ \ \ \ \ \ \ \ \ \ \ \ \ \ \ \ \ \ \ \ \ \ \ \ \ \ \ \ \ \ \ \ \ \ \ \ \ \ \ \nonumber
\end{eqnarray}

Therefore for every n$\geq$2
\begin{eqnarray}
\lefteqn {\rm d_n\leq2^{n{-}1}\ (p_1^{\frac{1}{{p_1}{-}1}}p_2^{\frac{1}{{p_2}{-}1}}. . .p_j^{\frac{1}{{p_j}{-}1}})^{n{-}2}\  ({2n{-}3})^{j/2} \ \ \ }  \ \ \ \ \ \ \ \ \ \ \ \ \ \ \ \ \ \ \ \ \ \ \ \ \ \ \ \ \ \ \ \ \ \ \ \ \ \ \ \ \ \ \ \ \ \ \ \ \ \ \ \ \ \ \ \ \ \ \ \ \ \ \ \ \ \ \ \ \ \ \ \ \ \nonumber
\end{eqnarray}
% therefore dn+1 ≤ 2n for n ≥ 1; of course d1 = 1 = 20.
\end{proof}

We denote by ${\rm D_n}$ the reduced denominator of ${\rm x_n}$. We have:
\begin{eqnarray}
\lefteqn {\rm D_n\geq{\frac{a_{n{-}1}}{d_n}}\geq\frac{\sqrt{(n{-}1)!}}{2^{n{-}1}(p_1^{\frac{1}{{p_1}{-}1}}p_2^{\frac{1}{{p_2}{-}1}}. . .p_j^{\frac{1}{{p_j}{-}1}})^{n{-}2}({2n{-}3})^{j/2}}\ \ \ \forall n\geq2} \ \ \ \ \ \ \ \ \ \ \ \ \ \ \ \ \ \ \ \ \ \ \ \ \ \ \ \ \ \ \ \ \ \ \ \ \ \ \ \ \ \ \ \ \ \ \ \ \ \ \ \ \ \ \ \ \ \ \ \ \ \ \ \ \ \ \ \ \ \ \ \ \ 
\end{eqnarray}

For simplicity we denote the rightmost side of (9) by ${\rm E(n)}$ as follows:
\begin{eqnarray}
\lefteqn {\rm E(n){=}\frac{\sqrt{(n{-}1)!}}{2^{n{-}1}(p_1^{\frac{1}{{p_1}{-}1}}p_2^{\frac{1}{{p_2}{-}1}}. . .p_j^{\frac{1}{{p_j}{-}1}})^{n{-}2}({2n{-}3})^{j/2}}\ \ \ \forall n\geq2} \ \ \ \ \ \ \ \ \ \ \ \ \ \ \ \ \ \ \ \ \ \ \ \ \ \ \ \ \ \ \ \ \ \ \ \ \ \ \ \ \ \ \ \ \ \ \ \ \ \ \ \ \ \ \ \ \ \ \ \ \ \ \ \ \ \ \ \ \ \ \ \ \ 
\end{eqnarray}

By inspecting, we have
\begin{eqnarray}
\lefteqn {\rm \lim_{n \to \infty}E(n){=}\lim_{n \to \infty} \frac{\sqrt{(n{-}1)!}}{2^{n{-}1}(p_1^{\frac{1}{{p_1}{-}1}}p_2^{\frac{1}{{p_2}{-}1}}. . .p_j^{\frac{1}{{p_j}{-}1}})^{n{-}2}({2n{-}3})^{j/2}}\ {=}\ \infty} \ \ \ \ \ \ \ \ \ \ \ \ \ \ \ \ \ \ \ \ \ \ \ \ \ \ \ \ \ \ \ \ \ \ \ \ \ \ \ \ \ \ \ \ \ \ \ \ \ \ \ \ \ \ \ \ \ \ \ \ \ \ \ \ \ \ \ \ \ \ \ \ \ 
\end{eqnarray}

It follows (11) that
\begin{eqnarray}
\lefteqn {\rm \exists h({\in}{\mathbb{N}})\ \ D_n{\geq}2\ \ \forall n{\geq}h} \ \ \ \ \ \ \ \ \ \ \ \ \ \ \ \ \ \ \ \ \ \ \ \ \ \ \ \ \ \ \ \ \ \ \ \ \ \ \ \ \ \ \ \ \ \ \ \ \ \ \ \ \ \ \ \ \ \ \ \ \ \ \ \ \ \ \ \ \ \ \ \ \ 
\end{eqnarray}

This means that if n${\geq}$h then ${\rm x_n}$ is never an integer.\\

At <First Solution> we conclude that the integral values of the sequence (${\rm x_n}$) are only ${\rm x_0}$ and ${\rm x_1}$. 
So if n${\geq}$2, then ${\rm D_n{\geq}2}$. 

If (12) hold for the case h=2, it is perfect. But we cannot at once judge that ${\rm E(n){\geq}}$2 for the case n${\geq}$2.\\

For example on Problem 1, c=1 and the inequality below holds.
\begin{eqnarray}
\lefteqn {\rm E(n){=}\frac{\sqrt{(n{-}1)!}}{2^{n{-}1}}\geq2\ \ \ \forall n\geq10} ,\ \ \ \ \ \ \ \ \ \ \ \ \ \ \ \ \ \ \ \ \ \ \ \ \ \ \ \ \ \ \ \ \ \ \ \ \ \ \ \ \ \ \ \ \ \ \ \ \ \ \ \ \ \ \ \ \ \ \ \ \ \ \ \ \ \ \ \ \ \ \ \ \ 
\end{eqnarray}

It is evident that ${\rm x_n}$ is not integral for n${\geq}$10. But we are still left to inspect the value of ${\rm x_n}$ with 0${\leq}$n${\leq}$9, as follows; 

${\rm x_0{=}}$1, ${\rm x_1{=}}$1, ${\rm x_2{=}}$2, ${\rm x_3{=}}$2, ${\rm x_4{=}}$$\frac{5}{2}$, ${\rm x_5{=}}$$\frac{13}{5}$, ${\rm x_6{=}}$$\frac{38}{13}$, ${\rm x_7{=}}$$\frac{58}{19}$, ${\rm x_8{=}}$$\frac{191}{58}$ and ${\rm x_9{=}}$$\frac{655}{191}$.

Thus by inspecting the change of E(n) according to n and computing its values for small n's, it can be known that the integral values of the sequence (${\rm x_n}$) are ${\rm x_0}$, ${\rm x_1}$, ${\rm x_2}$ and ${\rm x_3}$.\\

Similarly on Problem 2, for c$\geq$2, if h on (12) is found, then we are to get the solution only by computing ${\rm x_0}$, ${\rm x_1}$. . .${\rm x_{h{-}1}}$. 

But for now, it is far from easy to find h as a function of c. Especially for the cases that j is great or ${\rm p_1^{\frac{1}{{p_1}{-}1}}p_2^{\frac{1}{{p_2}{-}1}}. . .p_j^{\frac{1}{{p_j}{-}1}}}$ is great, it is hard to compute h. \\

Owing to (4), we can obtain another expression as follows;
\begin{eqnarray}
\lefteqn {\rm D_n\geq{\frac{a_{n{-}1}}{d_n}}{>}\frac{c^{n{-}1}}{2^{n{-}1}(p_1^{\frac{1}{{p_1}{-}1}}p_2^{\frac{1}{{p_2}{-}1}}. . .p_j^{\frac{1}{{p_j}{-}1}})^{n{-}2}({2n{-}3})^{j/2}}\ \ \ \forall n\geq2} \ \ \ \ \ \ \ \ \ \ \ \ \ \ \ \ \ \ \ \ \ \ \ \ \ \ \ \ \ \ \ \ \ \ \ \ \ \ \ \ \ \ \ \ \ \ \ \ \ \ \ \ \ \ \ \ \ \ \ \ \ \ \ \ \ \ \ \ \ \ \ \ \ \ \ \ \ \ 
\end{eqnarray}

Now suppose it holds ${\rm p_1}$<${\rm p_2}$<. . .<${\rm p_j}$. 

Then we consider each case according to the smallest prime factor of c.\\

\noindent <For the case of 2>\\

In this case c $\geq$ 2${\rm p_1}$${\rm p_2}$. . .${\rm p_j}$, so 
\begin{eqnarray}
\lefteqn {\rm \frac{c^{n{-}1}}{2^{n{-}1}(p_1^{\frac{1}{{p_1}{-}1}}p_2^{\frac{1}{{p_2}{-}1}}. . .p_j^{\frac{1}{{p_j}{-}1}})^{n{-}2}({2n{-}3})^{j/2}}{=}\frac{c^{n{-}1}}{2^{\frac{n}{2}}(\sqrt{2}p_1^{\frac{1}{{p_1}{-}1}}p_2^{\frac{1}{{p_2}{-}1}}. . .p_j^{\frac{1}{{p_j}{-}1}})^{n{-}2}({2n{-}3})^{j/2}}} \ \ \ \ \ \ \ \ \ \ \ \ \ \ \ \ \ \ \ \ \ \ \ \ \ \ \ \ \ \ \ \ \ \ \ \ \ \ \ \ \ \ \ \ \ \ \ \ \ \ \ \ \ \ \ \ \ \ \ \ \ \ \ \ \ \ \ \ \ \ \ \ \ \ \ \ \ \ \ \ \ \ \ \ \ \ \ \ \ \ \ \ \ \ \ \ \ \ \ \ \ \ \ \ \ \ \ \ \ \ \ \ \ \  \nonumber \\
\lefteqn {\rm \geq\frac{c^{n{-}1}}{2^{\frac{n}{2}}c^{\frac{n}{2}{-}1}({2n{-}3})^{j/2}}{=}{\left(\frac{c}{2}\right)}^{\frac{n}{2}}\frac{1}{{(2n{-}3})^{j/2}}} \ \ \ \ \ \ \ \ \ \ \ \ \ \ \ \ \ \ \ \ \ \ \ \ \ \ \ \ \ \ \ \ \ \ \ \ \ \ \ \ \ \ \ \ \ \ \ \ \ \ \ \ \ \ \ \ \ \ \ \ \ \ \ \ \ \ \ \ \ \ \ \ \ \ \ \ \ \ \ \ \ \ \ \ \ \ \ \ \ \ \ \ \ \ \ \ \ \ \ \ \ \ \ \ \ \ \ \ \ \ \ 
\end{eqnarray}

We denote the rightmost side of (15) by ${\rm E_2(n)}$ as follows:
\begin{eqnarray}
\lefteqn {\rm E_2(n){:=}\left(\frac{c}{2}\right)^{\frac{n}{2}}\frac{1}{{(2n{-}3})^{j/2}}} \ \ \ \ \ \ \ \ \ \ \ \ \ \ \ \ \ \ \ \ \ \ \ \ \ \ \ \ \ \ \ \ \ \ \ \ \ \ \ \ \ \ \ \ \ \ \ \ \ \ \ \ \ \ \ \ \ \ \ \ \ \ \ \ \ \ \ \ \ \ \ \ \ \ \ \ \ \ \ \ \ \ \ \ \ \ \ \ \ \ \ \ \ \ \ \ \ \ \ \ \ \ \ \ \ \ \ \ \ \ \ \nonumber 
\end{eqnarray}

Then ${\rm E_2(2)=\left(\frac{c}{2}\right)^1\frac{1}{1^1}=\frac{c}{2}\geq1}$ and
\begin{eqnarray}
\lefteqn {\rm \frac{E_2(n{+}1)}{E_2(n)}=\left(\frac{c}{2}\right)^{\frac{1}{2}}\left(\frac{2n{-}3}{2n{-}1}\right)^{j/2}\geq (p_1p_2. . .p_j)^\frac{1}{2}\left(\frac{1}{3}\right)^{j/2}\geq1} \ \ \ \ \ \ \ \ \ \ \ \ \ \ \ \ \ \ \ \ \ \ \ \ \ \ \ \ \ \ \ \ \ \ \ \ \ \ \ \ \ \ \ \ \ \ \ \ \ \ \ \ \ \ \ \ \ \ \ \ \ \ \ \ \ \ \ \ \ \ \ \ \ \ \ \ \ \ \ \ \ \ \ \ \ \ \ \ \ \ \ \ \ \ \ \ \ \ \ \ \ \ \ \ \ \ \ \ \ \ \ \nonumber
\end{eqnarray}

Therefore ${\rm E_2(n)\geq1}$ for $\forall$n$\geq$2, which results in ${\rm D_n{>}1}$ where n$\geq$2. \\

\noindent <For the case of over or equal to 5>\\

In this case ${\rm p_1{\geq}5}$, so
\begin{eqnarray}
\lefteqn {\rm \frac{c^{n{-}1}}{2^{n{-}1}(p_1^{\frac{1}{{p_1}{-}1}}p_2^{\frac{1}{{p_2}{-}1}}. . .p_j^{\frac{1}{{p_j}{-}1}})^{n{-}2}({2n{-}3})^{j/2}}\geq\frac{c^{n{-}1}}{2^{n{-}1}(p_1^{\frac{1}{4}}p_2^{\frac{1}{4}}. . .p_j^{\frac{1}{4}})^{n{-}2}({2n{-}3})^{j/2}}} \ \ \ \ \ \ \ \ \ \ \ \ \ \ \ \ \ \ \ \ \ \ \ \ \ \ \ \ \ \ \ \ \ \ \ \ \ \ \ \ \ \ \ \ \ \ \ \ \ \ \ \ \ \ \ \ \ \ \ \ \ \ \ \ \ \ \ \ \ \ \ \ \ \ \ \ \ \ \ \ \ \ \ \ \ \ \ \ \ \ \ \ \ \ \ \ \ \ \ \ \ \ \ \ \ \ \ \ \ \ \ \ \ \  \nonumber \\
\lefteqn {\rm \geq\frac{c^{n{-}1}}{2^{n{-}1}(c^{\frac{1}{4}})^{n{-}2}({2n{-}3})^{j/2}}{=}\frac{c^{\frac{3}{4}n{-}\frac{1}{2}}}{2^{n{-}1}({2n{-}3})^{j/2}}} \ \ \ \ \ \ \ \ \ \ \ \ \ \ \ \ \ \ \ \ \ \ \ \ \ \ \ \ \ \ \ \ \ \ \ \ \ \ \ \ \ \ \ \ \ \ \ \ \ \ \ \ \ \ \ \ \ \ \ \ \ \ \ \ \ \ \ \ \ \ \ \ \ \ \ \ \ \ \ \ \ \ \ \ \ \ \ \ \ \ \ \ \ \ \ \ \ \ \ \ \ \ \ \ \ \ \ \ \ \ \ \ \ \  
\end{eqnarray}

We denote the rightmost side of (16) by ${\rm E_3(n)}$ as follows:
\begin{eqnarray}
\lefteqn {\rm E_3(n){:=}\frac{c^{\frac{3}{4}n{-}\frac{1}{2}}}{2^{n{-}1}({2n{-}3})^{j/2}}} \ \ \ \ \ \ \ \ \ \ \ \ \ \ \ \ \ \ \ \ \ \ \ \ \ \ \ \ \ \ \ \ \ \ \ \ \ \ \ \ \ \ \ \ \ \ \ \ \ \ \ \ \ \ \ \ \ \ \ \ \ \ \ \ \ \ \ \ \ \ \ \ \ \ \ \ \ \ \ \ \ \ \ \ \ \ \ \ \ \ \ \ \ \ \ \ \ \ \ \ \ \ \ \ \ \ \ \ \ \ \ \nonumber 
\end{eqnarray}

Then ${\rm E_3(2)=\frac{c}{2}>1}$, and
\begin{eqnarray}
\lefteqn {\rm E_3(3)=\frac{c^{7/4}}{4(3^{j/2})}\geq \frac{(p_1p_2. . .p_j)^\frac{7}{4}}{2({\sqrt3})^j}>1} \ \ \ \ \ \ \ \ \ \ \ \ \ \ \ \ \ \ \ \ \ \ \ \ \ \ \ \ \ \ \ \ \ \ \ \ \ \ \ \ \ \ \ \ \ \ \ \ \ \ \ \ \ \ \ \ \ \ \ \ \ \ \ \ \ \ \ \ \ \ \ \ \ \ \ \ \ \ \ \ \ \ \ \ \ \ \ \ \ \ \ \ \ \ \ \ \ \ \ \ \ \ \ \ \ \ \ \ \ \ \ \nonumber \\
\lefteqn {\rm \frac{E_3(n{+}1)}{E_3(n)}=\frac{c^{3/4}}{2}\left(\frac{2n{-}3}{2n{-}1}\right)^{j/2}{\geq}\frac{1}{2}(p_1p_2. . .p_j)^\frac{3}{4}\left(\frac{3}{5}\right)^{j/2}{>}1\ \ \ \forall n\geq3} \ \ \ \ \ \ \ \ \ \ \ \ \ \ \ \ \ \ \ \ \ \ \ \ \ \ \ \ \ \ \ \ \ \ \ \ \ \ \ \ \ \ \ \ \ \ \ \ \ \ \ \ \ \ \ \ \ \ \ \ \ \ \ \ \ \ \ \ \ \ \ \ \ \ \ \ \ \ \ \ \ \ \ \ \ \ \ \ \ \ \ \ \ \ \ \ \ \ \ \ \ \ \ \ \ \ \ \ \ \ \ \nonumber
\end{eqnarray}

Therefore ${\rm E_3(n)\geq1}$ for $\forall$n$\geq$2, which results in ${\rm D_n{>}1}$ where n$\geq$2. \\

\noindent <For the case of 3>\\

In this case ${\rm p_1}$=3 and the other factors are all greater than 3, if they exist.

We assume the exponent of 3 in the prime factorization of c is more than 1, then \\
${\rm c\geq3^2p_2. . .p_j\geq(3^{1/2}p_2^{1/4}. . .p_j^{1/4})^4\geq(p_1^{\frac{1}{{p_1}{-}1}}p_2^{\frac{1}{{p_2}{-}1}}. . .p_j^{\frac{1}{{p_j}{-}1}})^4}$, therefore
\begin{eqnarray}
\lefteqn {\rm \frac{c^{n{-}1}}{2^{n{-}1}(p_1^{\frac{1}{{p_1}{-}1}}p_2^{\frac{1}{{p_2}{-}1}}. . .p_j^{\frac{1}{{p_j}{-}1}})^{n{-}2}({2n{-}3})^{j/2}}\geq\frac{c^{n{-}1}}{2^{n{-}1}(p_1^{\frac{1}{4}}p_2^{\frac{1}{4}}. . .p_j^{\frac{1}{4}})^{n{-}2}({2n{-}3})^{j/2}}} \ \ \ \ \ \ \ \ \ \ \ \ \ \ \ \ \ \ \ \ \ \ \ \ \ \ \ \ \ \ \ \ \ \ \ \ \ \ \ \ \ \ \ \ \ \ \ \ \ \ \ \ \ \ \ \ \ \ \ \ \ \ \ \ \ \ \ \ \ \ \ \ \ \ \ \ \ \ \ \ \ \ \ \ \ \ \ \ \ \ \ \ \ \ \ \ \ \ \ \ \ \ \ \ \ \ \ \ \ \ \ \ \ \  \nonumber \\
\lefteqn {\rm \geq\frac{c^{n{-}1}}{2^{n{-}1}(c^{\frac{1}{4}})^{n{-}2}({2n{-}3})^{j/2}}{=}\frac{c^{\frac{3}{4}n{-}\frac{1}{2}}}{2^{n{-}1}({2n{-}3})^{j/2}}} \ \ \ \ \ \ \ \ \ \ \ \ \ \ \ \ \ \ \ \ \ \ \ \ \ \ \ \ \ \ \ \ \ \ \ \ \ \ \ \ \ \ \ \ \ \ \ \ \ \ \ \ \ \ \ \ \ \ \ \ \ \ \ \ \ \ \ \ \ \ \ \ \ \ \ \ \ \ \ \ \ \ \ \ \ \ \ \ \ \ \ \ \ \ \ \ \ \ \ \ \ \ \ \ \ \ \ \ \ \ \ \ \ \  
\end{eqnarray}

Now we use ${\rm E_3(n)}$ again, then similarly as (15) ${\rm E_3(2){>}1}$, ${\rm E_3(3){>}1}$ and ${\rm \frac{E_3(n{+}1)}{E_3(n)}{>}1\ \ \forall n\geq3}$, which results in ${\rm D_n{>}1}$ where n$\geq$2.

Next we inspect the case in which the exponent of 3 is 1, then still ${\rm p_1}$=3.

Besides we assume that j>1, so ${\rm p_2{\geq}3}$ and c${\geq}$15.
\begin{eqnarray}
\lefteqn {\rm \frac{c^{n{-}1}}{2^{n{-}1}\left(p_1^{\frac{1}{{p_1}{-}1}}p_2^{\frac{1}{{p_2}{-}1}}. . .p_j^{\frac{1}{{p_j}{-}1}}\right)^{n{-}2}({2n{-}3})^{j/2}}{=}\frac{c^{n{-}1}}{2^{n{-}1}\left(\sqrt{3}p_2^{\frac{1}{{p_2}{-}1}}. . .p_j^{\frac{1}{{p_j}{-}1}}\right)^{n{-}2}({2n{-}3})^{j/2}}} \ \ \ \ \ \ \ \ \ \ \ \ \ \ \ \ \ \ \ \ \ \ \ \ \ \ \ \ \ \ \ \ \ \ \ \ \ \ \ \ \ \ \ \ \ \ \ \ \ \ \ \ \ \ \ \ \ \ \ \ \ \ \ \ \ \ \ \ \ \ \ \ \ \ \ \ \ \ \ \ \ \ \ \ \ \ \ \ \ \ \ \ \ \ \ \ \ \ \ \ \ \ \ \ \ \ \ \ \ \ \ \ \ \ \ \ \  \nonumber \\
\lefteqn {\rm \geq\frac{c^{n{-}1}}{2^{n{-}1}\left({\left(3^{\frac{2.5}{2}}p_2^{\frac{2.5}{{p_2}{-}1}}. . .p_j^{\frac{2.5}{{p_j}{-}1}}\right)}^{\frac{1}{2.5}}\right)^{n{-}2}({2n{-}3})^{j/2}}} \ \ \ \ \ \ \ \ \ \ \ \ \ \ \ \ \ \ \ \ \ \ \ \ \ \ \ \ \ \ \ \ \ \ \ \ \ \ \ \ \ \ \ \ \ \ \ \ \ \ \ \ \ \ \ \ \ \ \ \ \ \ \ \ \ \ \ \ \ \ \ \ \ \ \ \ \ \ \ \ \ \ \ \ \ \ \ \ \ \ \ \ \ \ \ \ \ \ \ \ \ \ \ \ \ \ \ \ \ \ \ \ \ \ \nonumber \\
\lefteqn {\rm \geq\frac{c^{n{-}1}}{2^{n{-}1}\left({\left(3^{1}p_2. . .p_j\right)}^{\frac{1}{2.5}}\right)^{n{-}2}({2n{-}3})^{j/2}}\geq\frac{c^{n{-}1}}{2^{n{-}1}\left({c}^{\frac{2}{5}}\right)^{n{-}2}({2n{-}3})^{j/2}}} \ \ \ \ \ \ \ \ \ \ \ \ \ \ \ \ \ \ \ \ \ \ \ \ \ \ \ \ \ \ \ \ \ \ \ \ \ \ \ \ \ \ \ \ \ \ \ \ \ \ \ \ \ \ \ \ \ \ \ \ \ \ \ \ \ \ \ \ \ \ \ \ \ \ \ \ \ \ \ \ \ \ \ \ \ \ \ \ \ \ \ \ \ \ \ \ \ \ \ \ \ \ \ \ \ \ \ \ \ \ \ \ \ \  \nonumber \\
\lefteqn {\rm {=}\frac{c^{\frac{3}{5}n{-}\frac{1}{5}}}{2^{n{-}1}({2n{-}3})^{j/2}}} \ \ \ \ \ \ \ \ \ \ \ \ \ \ \ \ \ \ \ \ \ \ \ \ \ \ \ \ \ \ \ \ \ \ \ \ \ \ \ \ \ \ \ \ \ \ \ \ \ \ \ \ \ \ \ \ \ \ \ \ \ \ \ \ \ \ \ \ \ \ \ \ \ \ \ \ \ \ \ \ \ \ \ \ \ \ \ \ \ \ \ \ \ \ \ \ \ \ \ \ \ \ \ \ \ \ \ \ \ \ \ \ \ \ 
\end{eqnarray}

We denote the rightmost side of (18) by ${\rm E_4(n)}$ as follows:
\begin{eqnarray}
\lefteqn {\rm E_4(n){:=}\frac{c^{\frac{3}{5}n{-}\frac{1}{5}}}{2^{n{-}1}({2n{-}3})^{j/2}}} \ \ \ \ \ \ \ \ \ \ \ \ \ \ \ \ \ \ \ \ \ \ \ \ \ \ \ \ \ \ \ \ \ \ \ \ \ \ \ \ \ \ \ \ \ \ \ \ \ \ \ \ \ \ \ \ \ \ \ \ \ \ \ \ \ \ \ \ \ \ \ \ \ \ \ \ \ \ \ \ \ \ \ \ \ \ \ \ \ \ \ \ \ \ \ \ \ \ \ \ \ \ \ \ \ \ \ \ \ \ \ \nonumber 
\end{eqnarray}

Then ${\rm E_4(2)=\frac{c}{2}>1}$, and
\begin{eqnarray}
\lefteqn {\rm E_3(3)=\frac{c^{8/5}}{4(3^{j/2})}\geq \frac{(p_1p_2. . .p_j)^\frac{8}{5}}{4({\sqrt3})^j}>1} \ \ \ \ \ \ \ \ \ \ \ \ \ \ \ \ \ \ \ \ \ \ \ \ \ \ \ \ \ \ \ \ \ \ \ \ \ \ \ \ \ \ \ \ \ \ \ \ \ \ \ \ \ \ \ \ \ \ \ \ \ \ \ \ \ \ \ \ \ \ \ \ \ \ \ \ \ \ \ \ \ \ \ \ \ \ \ \ \ \ \ \ \ \ \ \ \ \ \ \ \ \ \ \ \ \ \ \ \ \ \ \nonumber \\
\lefteqn {\rm \frac{E_3(n{+}1)}{E_3(n)}=\frac{c^{3/5}}{2}\left(\frac{2n{-}3}{2n{-}1}\right)^{j/2}{\geq}\frac{c^{3/5}}{2}\left(\frac{3}{5}\right)^{j/2}{>}1\ \ \ \forall n{\geq}3} \ \ \ \ \ \ \ \ \ \ \ \ \ \ \ \ \ \ \ \ \ \ \ \ \ \ \ \ \ \ \ \ \ \ \ \ \ \ \ \ \ \ \ \ \ \ \ \ \ \ \ \ \ \ \ \ \ \ \ \ \ \ \ \ \ \ \ \ \ \ \ \ \ \ \ \ \ \ \ \ \ \ \ \ \ \ \ \ \ \ \ \ \ \ \ \ \ \ \ \ \ \ \ \ \ \ \ \ \ \ \ \nonumber
\end{eqnarray}

Therefore ${\rm E_4(n)\geq1}$ for $\forall$n$\geq$2, which results in ${\rm D_n{>}1}$ where n$\geq$2. 

At last we are to inspect the case in which j=1, that is to say, c=3. In this case according to (10), 
\begin{eqnarray}
\lefteqn {\rm E(n){=}\frac{\sqrt{(n{-}1)!}}{2^{n{-}1}(\sqrt{3})^{n{-}2}({2n{-}3})^{1/2}}\ \ \ \forall n{\geq}2} \ \ \ \ \ \ \ \ \ \ \ \ \ \ \ \ \ \ \ \ \ \ \ \ \ \ \ \ \ \ \ \ \ \ \ \ \ \ \ \ \ \ \ \ \ \ \ \ \ \ \ \ \ \ \ \ \ \ \ \ \ \ \ \ \ \ \ \ \ \ \ \ \ \ \ \ \ \ \ \ \ \ \ \\
\lefteqn {\rm \frac{E(n{+}1)}{E(n)}{=}\frac{\sqrt{n}({2n{-}3})^{1/2}}{2\sqrt{3}({2n{-}1})^{1/2}}\ \ \ \forall n{\geq}2} \ \ \ \ \ \ \ \ \ \ \ \ \ \ \ \ \ \ \ \ \ \ \ \ \ \ \ \ \ \ \ \ \ \ \ \ \ \ \ \ \ \ \ \ \ \ \ \ \ \ \ \ \ \ \ \ \ \ \ \ \ \ \ \ \ \ \ \ \ \ \ \ \ \ \ \ \ \ \ \ \ \ \ 
\end{eqnarray}

Therefore it follows that 
\begin{eqnarray}
\lefteqn {\rm (E(30))^2{=}\frac{\sqrt{29!}}{4^{29}\ 3^{28}\ 57}{\geq}\frac{1.65\times10^{17}}{3.90\times10^{17}}{\geq}0.423} \ \ \ \ \ \ \ \ \ \ \ \ \ \ \ \ \ \ \ \ \ \ \ \ \ \ \ \ \ \ \ \ \ \ \ \ \ \ \ \ \ \ \ \ \ \ \ \ \ \ \ \ \ \ \ \ \ \ \ \ \ \ \ \ \ \ \ \ \ \ \ \ \ \ \ \ \ \ \ \ \ \ \ \\
\lefteqn {\rm \left(\frac{E(31)}{E(30)}\right)^2{=}\frac{30\times57}{12\times59}{\geq}2.415} \ \ \ \ \ \ \ \ \ \ \ \ \ \ \ \ \ \ \ \ \ \ \ \ \ \ \ \ \ \ \ \ \ \ \ \ \ \ \ \ \ \ \ \ \ \ \ \ \ \ \ \ \ \ \ \ \ \ \ \ \ \ \ \ \ \ \ \ \ \ \ \ \ \ \ \ \ \ \ \ \ \ \ 
\end{eqnarray}

It follows that ${\rm (E(31))^2{>}1}$, so ${\rm E(31){>}1}$. And on (20) if n>30, then ${\rm \frac{E(n{+}1)}{E(n)}{>}1}$, which results in ${\rm E(n){>}1}$ where ${\rm n{\geq}31}$. Sufficiently it holds that ${\rm D_n{>}1}$ where ${\rm n{\geq}31}$.

We have only to inspect ${\rm D_n}$ where ${\rm n{\leq}30}$.

With (4), for ${\rm 1{\leq}\forall n{\leq}15}$, ${\rm a_n}$ and ${\rm d_n{=}gcd(a_n,a_{n{-}1})}$ can be obtained as follows.
\begin{figure}[htbp]
\begin{center}
\includegraphics[width=5.5cm]{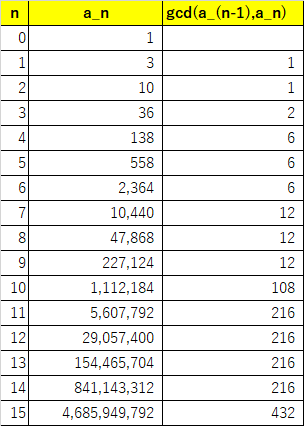}
% \caption{Author}
\end{center}
\end{figure}
\newpage
${\rm D_n}$ is evidently much more than 1 as above. 

Next for ${\rm 16{\leq}\forall n{\leq}30}$, it holds that
\begin{eqnarray}
\lefteqn {\rm D_n\geq{\frac{a_{n{-}1}}{d_n}}{>}\frac{a_{15}\times3^{n{-}{16}}}{2^{n{-}1}(\sqrt3)^{n{-}2}{\sqrt{2n{-}3}}}{=}\frac{a_{15}}{2^{15}(\sqrt3)^{14}{\sqrt{2n{-}3}}}\left(\frac{\sqrt3}{2}\right)^{n{-}16}} \ \ \ \ \ \ \ \ \ \ \ \ \ \ \ \ \ \ \ \ \ \ \ \ \ \ \ \ \ \ \ \ \ \ \ \ \ \ \ \ \ \ \ \ \ \ \ \ \ \ \ \ \ \ \ \ \ \ \ \ \ \ \ \ \ \ \ \ \ \ \ \ \ \ \ \ \ \ \ \ \ \ \ \ \ \ \ \ \ \ \ \ \ \ \ \ \ \ \ \ \ \ \ \ \nonumber \\
\lefteqn {\rm {>}\frac{a_{15}}{2^{15}(\sqrt3)^{14}{\sqrt{57}}}\left(\frac{\sqrt3}{2}\right)^{14}{=}\frac{a_{15}}{2^{29}{\sqrt{57}}}} \ \ \ \ \ \ \ \ \ \ \ \ \ \ \ \ \ \ \ \ \ \ \ \ \ \ \ \ \ \ \ \ \ \ \ \ \ \ \ \ \ \ \ \ \ \ \ \ \ \ \ \ \ \ \ \ \ \ \ \ \ \ \ \ \ \ \ \ \ \ \ \ \ \ \ \ \ \ \ \ \ \ \ \ \ \ \ \ \ \ \ \ \ \ \ \ \ \ \ \ \nonumber \\
\lefteqn {\rm {>}\frac{4,685,949,792}{536,870,912{\times}{7.55}}{>}\frac{4.68\times10^9}{4.06\times10^9}} \ \ \ \ \ \ \ \ \ \ \ \ \ \ \ \ \ \ \ \ \ \ \ \ \ \ \ \ \ \ \ \ \ \ \ \ \ \ \ \ \ \ \ \ \ \ \ \ \ \ \ \ \ \ \ \ \ \ \ \ \ \ \ \ \ \ \ \ \ \ \ \ \ \ \ \ \ \ \ \ \ \ \ \ \ \ \ \ \ \ \ \ \ \ \ \ \ \ \ \ \nonumber \\
\lefteqn {\rm {>}1.15} \ \ \ \ \ \ \ \ \ \ \ \ \ \ \ \ \ \ \ \ \ \ \ \ \ \ \ \ \ \ \ \ \ \ \ \ \ \ \ \ \ \ \ \ \ \ \ \ \ \ \ \ \ \ \ \ \ \ \ \ \ \ \ \ \ \ \ \ \ \ \ \ \ \ \ \ \ \ \ \ \ \ \ \ \ \ \ \ \ \ \ \ \ \ \ \ \ \ \ \ \nonumber
\end{eqnarray}

Therefore if c=3, then ${\rm D_n{>}1}$ $\forall$n$\geq$2. \\

We conclude that the only integral values of the sequence are ${\rm x_0}$, ${\rm x_1}$, similarly as seen in <First Solution>.\\
% We now prove the following lemma.
% {\bf Lemma 6.} \\\\
% \lim_{x \to \infty} f(x)
% e (an)n∈N, namely
% \begin{proof}
% \quad\par
% So (40) follows easily.\\
% \end{proof}
 %
% \\
% \begin{theorem} %Theorem 12
% \begin{eqnarray}
% && {\rm \lim_{n \to \infty} \dfrac{\left|C(n)\right|}{\left|S(n)\right|}}{=}0
% \end{eqnarray}
% \end{theorem} %Theorem 12
\\\\
% \noindent 4. Discussion and conclusion\\
%
%At <First Research> we show the solution successfully by the reference of [1].
%
%But at <Second Research> we do not, where we only see that there exists an positive integer h, such that ${\rm x_n}$ is never an integer if n$\geq$h, which is far from solving Problem 2. 
%
%In the future we want to solve it by the method at <Second Research>.
% By diagonalization a 5$\times$5 matrix, we have succeeding in expressing $\left|{\rm D^*(k)} \right|$ by the polynomials above.\\
\newpage
% \begin{table}[h]
% \caption{}
% \centering
% \small
% \begin{tabular}{|c|l|r|r|r|r|r|r|r|r|r|r|}\hline
%n&(=$\left|{\rm E(n)}\right|$) & 1 & 2 & 3 & 4 & 5 & 6 & 7 & 8 & 9 & 10\\ \hline
%$\left|{\rm F(n)}\right|$ & & 1 & 3 & 7 & 15 & 31 & 63 & 127 & 255 & 511 & 1023 \\ \hline
%  & $\left|{\rm C(n)}\right|$ & 0 & 3 & 20 & 95 & 384 & 1557 & 6101 & 23519 & 89924 & 342779 \\ \cline{2-12}
%  & $\left|{\rm C(n)}\right|$/$\left|{\rm S(n)}\right|$ & - & 1.00 & 0.952 & 0.905 & 0.826 & 0.797 & 0.763 & 0.726 & 0.690 &0.656 \\ \cline{2-12}
% $\left|{\rm P(n,1)}\right|$ & $\left|{\rm P_C(n,1)}\right|$ & 0 & 2 & 3 & 4 & 5 & 6 & 7 & 8&9&10\\ \cline{2-12}
%Then it can be easily understood that what we have described on this paper holds exactly the same.
% \\\\
\centerline{References}

\noindent [1] Roberto Dvornicich, Francesco Veneziano, and Umberto Zannier, \\6, arvix, available at https://arxiv.org/abs/1403.3470\\
\noindent [2] Romanian Mathematical Society, RMC 2010, Paralela 45 Publishing House, Bucharest (2010).\\
\end{document}